\documentclass[12pt]{article}

\usepackage{amsfonts,amsbsy,graphicx}

\def\blacksquare{\vrule height .9ex width .8ex depth -.1ex}

\newcommand{\ih}{\'{\i}}
\newcommand{\eh}{\hspace{.06in}}
\newcommand{\D}{{\cal{D}}}
\newcommand{\G}{{\cal{G}}}
\newcommand{\K}{{\cal{K}}}
\newcommand{\M}{{\cal{M}}}
\newcommand{\C}{\mathbb{C}}
\newcommand{\R}{\mathbb{R}}
\newcommand{\m}{\frac{1}{2}}
\newcommand{\af}{\alpha}
\newcommand{\de}{\delta}

\newcommand{\ta}{\tan\af}
\newcommand{\taa}{\tan^2\af}
\newcommand{\cta}{\cot\af}
\newcommand{\ctaa}{\cot^2\af}
\newcommand{\lessim}{\lower1.0pt\hbox{${_{\ds<}}\atop^{^{\ds\sim}}$}}
\newcommand{\gtrsim}{\lower1.0pt\hbox{${_{\ds>}}\atop^{^{\ds\sim}}$}}
\newcommand{\ds}{\displaystyle}
\newcommand{\ovl}{\overline}
\newcommand{\BE}{\begin{equation}}
\newcommand{\EE}{\end{equation}}
\newcommand{\fns}{\footnotesize}
\newcommand{\Lim}[1]{\lower5pt\hbox{${{\ds\lim}\atop ^{#1}}$}}

\begin{document}
\centerline{\bf THE DOUBLY PERIODIC SCHERK-COSTA SURFACES}
\bigskip

\centerline{K{\fns ELLY} L{\fns \"UBECK} \& V{\fns AL\'ERIO} R{\fns AMOS} B{\fns ATISTA}}
\bigskip
\begin{abstract}
We present a new family of embedded doubly periodic minimal surfaces, of which the symmetry group does not coincide with any other example known before.
\end{abstract}
\ \\
{\bf 1. Introduction}
\\

In the Euclidean three-space, among all complete embedded minimal surfaces known to date, on the one hand the doubly periodic class still remains less numerous in examples. On the other hand, the richness in the triply periodic class counts to a great extent on the {\it Conjugate Plateau Construction}, a powerful tool but not always applicable (see [RB2]) or extendable to infinite frames (see [JS]). Recently, the non-periodic class became also very rich due to important works like [Kp] and [T3]. Through [T1-2] and [Web2] the same happened to the singly periodic, which was already quite numerous with tree possible kinds of ends: planar, Scherk or helicoidal. 
\\

In the 20$^{\rm th}$ century, the known examples were obtained thanks to their high order symmetry groups, a resource already used up nowadays. Therefore, potentially new examples normally lack in symmetries, which makes it so hard to prove their existence. One might opt for keeping a high order symmetry group with an increase in the genus, but this leads to the same hurdle, namely too many {\it period problems}.
\\

Some modern constructions {\it do} close many periods at once, like in [Kp], [T1-3] and [Web2]. However, such methods are not portable outside the class and types of ends they direct to. For at most three period problems, however, it is still feasible to handle the {\it Weierstra\ss \ Data} together with adequate methods (see [BRB], [L], [LRB], [MRB] and [Web1]). By the way, in our paper we apply the {\it limit-method} described both in [L] and [LRB].      
\\

However, what is the purpose of constructing a new minimal surface? The following reason motivates this present work. Minimal surfaces model many structures, like crystals and co-polymers, but several symmetry groups are not yet represented by any them (see [H] and [LM]). As explained above, the doubly periodic class still lacks in examples, even after very rich works like [HKW], [K1], [W] and [PRT]. The purpose of this paper is then to present a new family of embedded doubly periodic minimal surfaces, of which the symmetry group does not coincide with any other example known before. For the converse, {\it there are} symmetry groups that admit more than one representative, even restricted to a certain conformal type (see [RB3]). Although not embedded, these examples easily hint at embedded ones.     
\\

Our surfaces are inspired in the examples called L$_{\rm b}$ in [RB4, p 482]. By taking the picture of L$_{\rm b}$ in that page, if one replaces the catenoidal ends by curves of reflectional symmetry parallel to $Ox_1x_2$, the resultant surface will then come out as in Figure 1. Its symmetry group is $D_{2d}$ or $\bar{4}2m$ in Sch\"onflies' or simplified Hermann-Maugin notation, respectively. The same procedure for C$_{\rm b}$ from [RB4, p 483] could also result in a new doubly periodic example. However, it would then have the same symmetry group as $M_{\pi/4,\pi/2,0}$ from [PRT]. Of course, they would differ in genus, but one still might go round it by {\it adding handles} to the latter, a widely applicable technique. That is why our paper is totally devoted to the example in Figure 1.
\\

We formally state our result in the following theorem: 
\begin{figure}[!ht]
 \centering
 \includegraphics[scale=0.7]{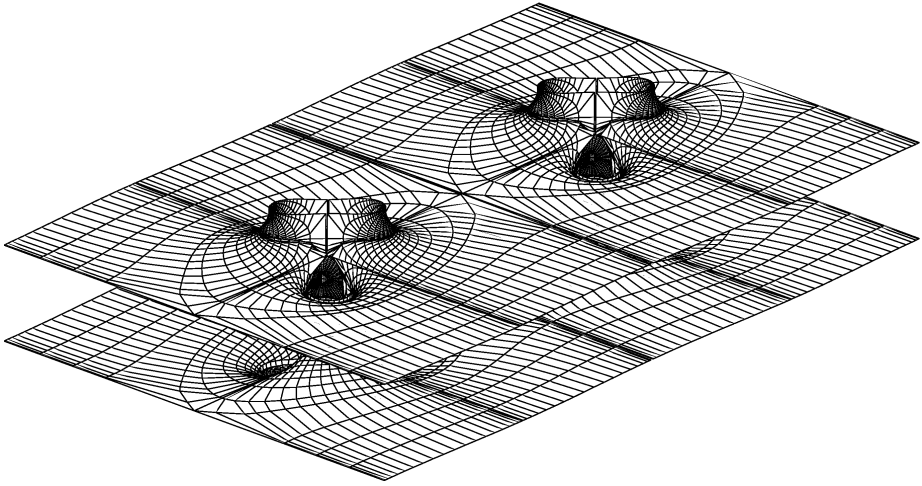}
\caption{A doubly periodic Scherk-Costa surface.}
\end{figure}
\ \\
{\bf Theorem 1.1}. \it There exists a one-parameter family of complete embedded doubly periodic minimal surfaces in $\R^3$ such that, for each member of this family the following holds:
\\
(a) The quotient by its translation group has genus three.
\\
(b) The surface is generated by a fundamental piece, which is a surface with boundary in $\R^3$. The fundamental piece has two Scherk-ends (modulo translation), and a symmetry group generated by 180$^\circ$-rotations around a straight line and 180$^\circ$-rotations around a straight segment. The segment crosses the line orthogonally and both determine a plane $\Pi$.
\\
(c) The boundary of the fundamental piece consists of two parallel lines in $\Pi$ and two planar closed curves of reflectional symmetry. The curves are parallel to but not contained in $\Pi$, and one is the image of the other under the symmetries of the fundamental piece. By successive 180$^\circ$-rotations around the lines of the boundary and reflections in the closed curves, one generates the doubly periodic surface.\rm
\\

This work refers to part of our doctoral theses [L] and [RB0], supported by CAPES (Coordena\c c\~ao de Aperfei\c coamento de Pessoal de N\ih vel Superior) and DAAD (Deutscher Akademischer Austausch Dienst), respectively. Professor Hermann Karcher, from the University of Bonn in Germany, was adviser of the second author, who thanks him for his dedication, which greatly helped in the realisation of this work. The second author was advisor of the first.   
\\
\\
{\bf 2. Preliminaries}
\\

In this section we state some well known theorems on minimal surfaces. For details, we refer the reader to [K2], [LMa], [N] and [O]. In this paper all surfaces are supposed to be regular.
\\
\\
{\bf Theorem 2.1.} (Weierstra\ss \ representation). \it Let $R$ be a Riemann surface, $g$ and $dh$ meromorphic function and 1-differential form on $R$, respectively, such that the zeros of $dh$ coincide with the poles and zeros of $g$. Consider the (possibly multi-valued) function $X:R\to\R^3$ given by
\[
   X(p):=Re\int^p(\phi_1,\phi_2,\phi_3),\eh\eh where\eh\eh
   (\phi_1,\phi_2,\phi_3):=\m(g^{-1}-g,ig^{-1}+ig,2)dh.
\]
Then $X$ is a conformal minimal immersion. Conversely, every conformal minimal immersion $X:R\to\R^3$ can be expressed as above for some meromorphic function $g$ and 1-form $dh$.\rm
\\
\\
{\bf Definition 2.1.} The pair $(g,dh)$ is the \it Weierstra\ss \ data \rm and $\phi_1$, $\phi_2$, $\phi_3$ are the \it Weierstra\ss \ forms \rm on $R$ of the minimal immersion $X:R\to X(R)=S\subset\R^3$.
\\
\\
{\bf Definition 2.2.} A complete, orientable minimal surface $S$ is \it algebraic \rm if it admits a Weierstra\ss \ representation such that $R=\ovl{R}\setminus\{p_1,\dots,p_r\}$, were $\ovl{R}$ is compact, and both $g$ and $dh$ extend meromorphically to $\ovl{R}$.
\\
\\
{\bf Definition 2.3.} An \it end \rm of $S$ is the image of a punctured neighbourhood $V_p$ of a point $p\in\{p_1,\dots,p_r\}$ such that $(\{p_1,\dots,p_r\}\setminus\{p\})\cap\ovl{V}_p=\emptyset$. The end is \it embedded \rm if this image is embedded for a sufficiently small neighbourhood of $p$.
\\
\\
{\bf Theorem 2.2.} \it Let $S$ be a complete minimal surface in $\R^3$. Then $S$ is algebraic if and only if it can be obtained from a piece $\tilde{S}$ of finite total curvature by applying a finitely generated translation group $G$ of \ $\R^3$.\rm
\\

From now on we consider only algebraic surfaces. The function $g$ is the stereographic projection of the Gau\ss \ map $N:R\to S^2$ of the minimal immersion $X$. This minimal immersion is well defined in $\R^3/G$, but allowed to be a multivalued function in $\R^3$. The function $g$ is a covering map of $\hat\C$ and the total curvature of $\tilde{S}$ is $-4\pi$deg$(g)$.
\\
\\
{\bf 3. The compact Riemann surfaces $\ovl{M}$ and the functions $z$}
\\

Denote by $M$ the surface represented in Figure 1, and let $\M$ be the quotient of $M$ by its translation group. A compactification of the Scherk ends of $\M$ will lead to a compact Riemann surface that we call $\ovl{M}$. The fundamental piece of $M$ is represented in Figure 2(a), together with some special points on it. The Scherk ends are $E_1$ and $E_2$. We have that $\ovl{M}$ is invariant under reflections in the closed bold curve, indicated in Figure 2(a). The images of $S,F,E_1$ and $E_2$ under this reflection will be called $S',F',E_1$ and $E_2$, respectively. 
\\
\input epsf
\begin{figure} [ht]
\centerline{
\epsfxsize 13cm
\epsfbox{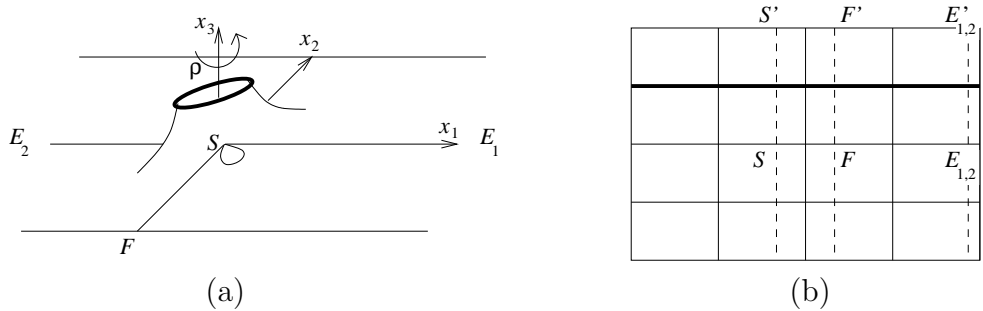}}
\hspace{1.2in}(a)\hspace{2.85in}(b)
\caption{(a) The fundamental piece of $M$; (b) the torus $T=\rho(\ovl{M})$.}
\end{figure}

Let us denote by $\rho$ the $180^\circ$-rotation around the axis $x_3$, indicated in Figure 2(a). It is easy to see that $\ovl{M}$ has genus 3 and $\rho$ has 4 fixed points, namely $S,T,S',F'$. Therefore, the Euler-Poincar\'e formula gives
\[
   \chi(\rho(\ovl{M}))=\frac{\chi(\ovl{M})}{2}+\frac{4}{2}=0.
\]     

Because of that, $\rho(\ovl{M})$ is a torus $T$. A horizontal reflectional symmetry of $\ovl{M}$ is induced by $\rho$ on $T$, and since this symmetry has two components, we conclude that $T$ is a rectangular torus, represented in Figure 2(b). Now we can choose an elliptic function $Z$ on $T$, define $z:=Z\circ\rho$ and then try to deduce Weierstra\ss \ data of $M$ on $\ovl{M}$ in terms of $z$. Consider Figure 3(a) and the points marked with a black square ($\blacksquare$) on it, which are the branch points of a certain meromorphic function $Z:T\to\C$ with deg($Z$)=2. Let us now take an angle $\af\in(0,\pi/4)$. As indicated in Figure 3(a), we choose $Z$ such that it takes the values $\pm i\ta$ and $\pm i\cta$ at its branch points. Up to a biholomorphism, such a function is unique, and $\af$ determines the rectangular torus. The torus is square if and only if $\af=\pi/8$.

\input epsf
\begin{figure} [ht]
\centerline{
\epsfxsize 13cm
\epsfbox{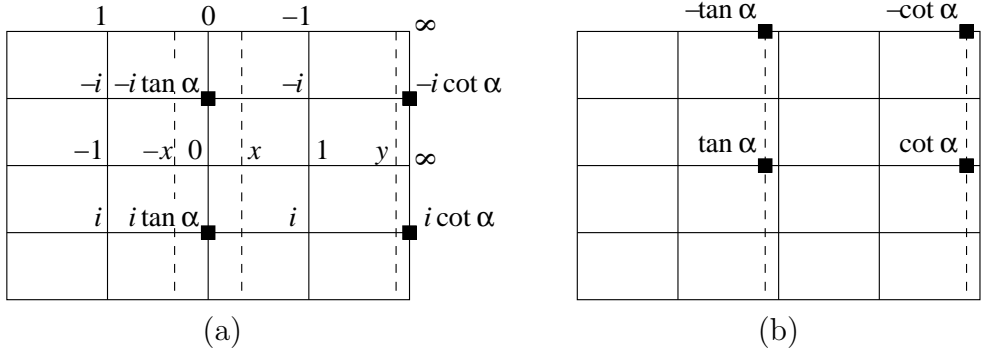}}
\hspace{1.2in}(a)\hspace{2.7in}(b)
\caption{(a) Important values of $Z$ on $T$; (b) the function $W$ on $T$.}
\end{figure}

The most important values of $Z$ are indicated in Figure 3(a). We have $Z(\rho(F))=-Z(\rho(S))=x$ for some $x\in(0,\infty)$, and $Z(\rho(E_{1,2}))=y$, where $y^{-1}\in(-x^{-1},x^{-1})$. This means, we include the possibility of $y$ to be $\infty$. Consequently, $Z(\rho(S'))=-Z(\rho(F'))=x$ and $Z(\rho(E_{1,2}))=-y$.
\\
 
In the next section, we shall write the function $g$ on $\ovl{M}$ in terms of $z:=Z\circ\rho$. However, this task will be simpler if we introduce another function $W:T\to\C$, of which the important values are presented in Figure 3(b). In fact, $W$ is a ``shift'' of $iZ$. We can write $W$ in terms of $Z$ and $Z^\prime$ by using an addition theorem for elliptic functions, or apply some {\it easier} arguments which will be explained in the next paragraph. Nevertheless, they will give us an explicit formula for $W^2$ instead of $W$. 
\\

Consider the following picture where $\xi$ is a pure imaginary value to be determined later:
\input epsf
\begin{figure} [ht]
\centerline{
\epsfxsize 13cm
\epsfbox{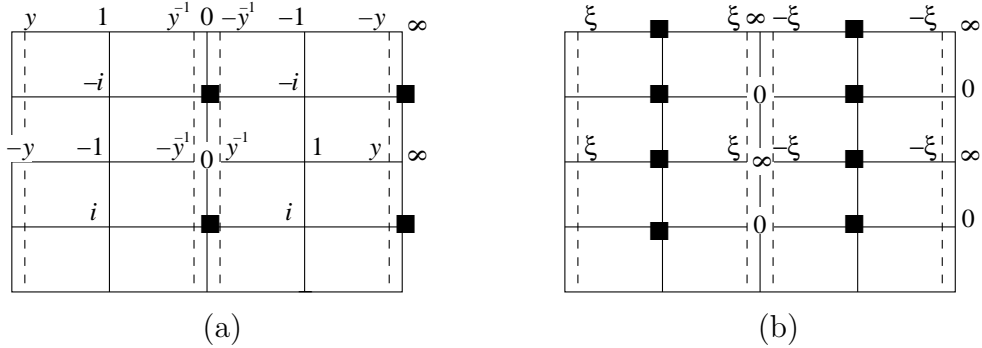}}
\hspace{1.2in}(a)\hspace{2.7in}(b)
\caption{(a) Poles and zeros of $Z$; (b) poles and zeros of $Z^\prime/Z$.}
\end{figure}

Let us define 
\[
   f:=\frac{Z^2-y^{-2}}{Z\biggl(\frac{\ds Z^\prime}{\ds Z}+\xi\biggl)}.
\]

From Figures 3(a) and 4 it is easy to see that, for a certain complex constant $c$, the equality $cf^2=\frac{\ds W^2-\taa}{\ds W^2-\ctaa}$ holds. Based on Figure 3(a), we can easily write down an algebraic equation of $T$ as follows:
\BE
   Z^{\prime 2}=-(Z^2+\taa)(Z^2+\ctaa).
\EE

Therefore, $\biggl(\frac{\ds Z^\prime}{\ds Z}\biggl)^2=-(Z^2+Z^{-2}+\taa+\ctaa)$ and consequently we fix $\xi=i(y^2+y^{-2}+\taa+\ctaa)^\m$. From Figure 3, it is not difficult to prove that $W=\infty$ implies $Z=\pm i\biggl(\frac{\ds 1+y^2\ctaa}{\ds y^2+\ctaa}\biggl)^\m$. Hence, 
\[
   c=f^{-2}|_{W=\infty}=\frac{\ctaa\cdot y^2}{[2+\ctaa\cdot(y^2+y^{-2})]^2}\cdot
   (\taa-\ctaa+|\xi|^2)^2.
\]

The explicit relation between $W^2$ and $Z,Z^\prime$ can be given as follows
\BE
   W^2=\ctaa+\frac{\ctaa-\taa}{cf^2-1},\eh\eh{\rm where}\eh\eh
   f=\frac{Z^2-y^{-2}}{Z^\prime+\xi Z}.
\EE
\ \\
{\bf 4. The Gau\ss \ map of $M$ and the function $g$ on $\ovl{M}$}
\\

Consider $z:=Z\circ\rho$ and $w:=W\circ\rho$. Therefore, $z$ and $w$ are meromorphic functions on $\ovl{M}$ and deg($z$)=deg($w$)=4. Based on Figures 1 and 2(a), one easily sees that the unitary normal vector on $M$ is expected to be vertical at $S,F,E_1$ and $E_2$. From now on, we are going to use some heuristic arguments: if the normal vector points downwards at $S$, it will then point upwards at $F,E_1$ and $E_2$. Consequently, $g(\{S,F',E_{1,2}\})=\{0\}$ and $g(\{S',F,E_{1,2}\})=\{\infty\}$. We do not expect the normal vector to be vertical at any other point of $\ovl{M}$, except at the ones just mentioned. Moreover, all the poles and zeros of $g$ are simple. Hence, deg($g$)=4.
\\
\input epsf
\begin{figure} [ht]
\centerline{
\epsfxsize 13cm
\epsfbox{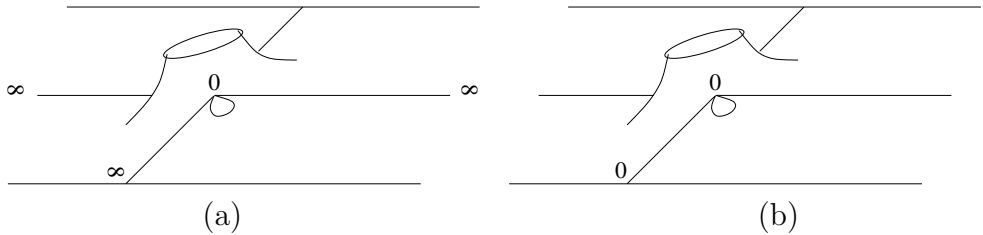}}
\hspace{1.2in}(a)\hspace{2.7in}(b)
\caption{(a) Poles and zeros of $g$; (b) poles and zeros of $dh$.}
\end{figure}

Based on Figures 3 and 5(a), one easily concludes that
\BE
   g^2=\frac{x+z}{x-z}\cdot\frac{\cta+w}{\cta-w}.
\EE  
   
Of course, a priori both sides of 3 are just proportional. However, since the unitary normal vector is expected to be horizontal on the closed bold curve in Figure 2(a), this must imply that $g$ is unitary there. Both $z$ and $w$ are pure imaginary on this curve. Therefore, the proportional constant at (3) must be unitary. Moreover, based on Figures 2(a) and 3, the first picture suggests that $g$ is real for $z\in(-x,x)$ and pure imaginary for $z\in\R\setminus(-x,x)$. Since $w\in[-\cta,\cta]$ whenever $z$ is real, we conclude that the unitary proportional constant is 1. Thus, (3) itself is already consistent with our analysis. From now on, we define $\ovl{M}$ as a member of the family of compact Riemann surfaces given by (3).
\\

By applying the Riemann-Hurwitz formula to (3), one obtains 
\[
   \frac{\ds 4(2-1)+8(2-1)}{\ds 2}-4+1=3.
\]
Thus, the genus of $\ovl{M}$ is three. Now we must verify that $\ovl{M}$ really has all the symmetries we supposed at the beginning, and $g$ really corresponds to the unitary normal vector on $M$. First of all, let us show that $w(\ovl{z})=\ovl{w}(z)$ and $w(-\ovl{z})=-\ovl{w}(z)$. From (1) we have $Z^\prime(\ovl{Z})=\pm\ovl{Z^\prime}(Z)$. But according to (1), $Z^\prime$ is pure imaginary where $Z$ is real. Hence $Z^\prime(\ovl{Z})=-\ovl{Z^\prime}(Z)$. Now we recall that $\xi=i|\xi|$, use (2) and get $f(\ovl{Z})=-\ovl{f}(Z)$ and $W^2(\ovl{Z})=\ovl{W}^2(Z)$. Hence $W(\ovl{Z})=\pm\ovl{W}(Z)$, but since $W\in[-\cta,\cta]$ whenever $Z$ is real, then $W(\ovl{Z})=\ovl{W}(Z)$. Now apply $z=Z\circ\rho$ and $w=W\circ\rho$ to these relations. 
\\

By recalling (1) again, in the case $Z\to-\ovl{Z}$ we do have two possibilities: either $Z^\prime(-\ovl{Z})=\ovl{Z^\prime}(Z)$ or $Z^\prime(-\ovl{Z})=-\ovl{Z^\prime}(Z)$. Nevertheless, our assumptions about the symmetries of $M$ {\it do not} imply that $\rho$ induces from $\ovl{M}$ the involution $(Z^\prime,Z)\to(-\ovl{Z^\prime},-\ovl{Z})$ on $T$. Hence, we only consider $(Z^\prime,Z)\to(\ovl{Z^\prime},-\ovl{Z})$. From (2) it follows that $W(-\ovl{Z})=\pm\ovl{W}(Z)$. Since $W$ must be pure imaginary for $Z\in i\cdot[-\cta,\cta]$, then $W(-\ovl{Z})=-\ovl{W}(Z)$ and consequently we get $w(-\ovl{z})=-\ovl{w}(z)$. 
\\

Now we can summarize our study of the symmetries of $\ovl{M}$ and the behaviour of $g$ in the following table:
\BE
\begin{tabular}{|c|c|c|c|}\hline
 ${\rm involution}             $&$ z{\rm -values}    $&$ g{\rm -values} $ \\ \hline\hline
 $(g,z)\to(\ovl{g},\ovl{z})    $&$ z=t,-x<t<x        $&$ g\in  \R    $ \\ \hline 
 $(g,z)\to(-\ovl{g},\ovl{z})   $&$ z=t,x<t<y         $&$ g\in i\R    $ \\ \hline
 $(g,z)\to(-\ovl{g},\ovl{z})   $&$ z=t,y<t<-x        $&$ g\in i\R    $ \\ \hline
 $(g,z)\to(1/\ovl{g},-\ovl{z}) $&$ z=it,\ta<t<\cta   $&$ |g|\equiv 1    $ \\ \hline
 $(g,z)\to(1/\ovl{g},-\ovl{z}) $&$ z=it,-\cta<t<-\ta $&$ |g|\equiv 1    $ \\ \hline
\end{tabular}
\EE
\\
\\
{\bf 5. The height differential $dh$ on $\ovl{M}$}
\\

Now we are going to write down an explicit formula for $dh$, which will take into account the regular points and types of ends we want the surface $M$ to have. Based on Figures 2(a) and 3(a), one sees that $S$ and $F$ correspond to regular points of $M$, at which the normal vector is vertical. The same is valid for $S'$ and $F'$. Therefore, $dh(\{S,F,S',F'\})=\{0\}$. Since $M$ has only Scherk-type ends, all in the $x_2$-direction, then $dh$ has no poles and is holomorphic on $\ovl{M}$. Moreover, since deg($dh$)=4, we conclude that all zeros of $dh$ are simple. They are represented in Figure 5(b).
\\

For convenience of the reader, we reproduce here the algebraic equation of $T$ already established in (1):
\BE
   Z^{\prime 2}=-(Z^2+\taa)(Z^2+\ctaa).
\EE

From Figure 5(b) one immediately verifies that
\BE
   dh=\frac{dz}{Z^\prime\circ\rho}.
\EE 

A priori, both sides of (6) are just proportional, but since $z$ is real on the straight lines of $\ovl{M}$, from (5) we get pure imaginary values for $Z^\prime\circ\rho$ on these lines. Therefore, the proportional constant at (6) must be real, and we choose it to be 1. This will imply that (6) is also consistent with $Re\{dh\}=0$ on the closed bold curve in Figure 2(a). There we have $z=it,\ta<|t|<\cta$, which leads to real values for $Z^\prime\circ\rho$.
\\

Analogously, we could also have defined the algebraic equation of $T$ by $W^{\prime 2}=(W^2-\taa)(W^2-\ctaa)$ and so $dh=dw/W^\prime\circ\rho$, namely
\BE
   dh=\frac{dw}{\sqrt{(w^2-\taa)(w^2 -\ctaa)}}.
\EE

Exactly at this point, we need to prove that $M$ really has the planar geodesics and straight lines of our initial assumption. This task is summarized in the following table:
\[
\begin{tabular}{|c|c|c|c|}\hline
 $ z=t,-x<t<x        $&$ g\in  \R    $&$ dh(\dot{z})\in i\R $ \\ \hline 
 $ z=t,x<t<y         $&$ g\in i\R    $&$ dh(\dot{z})\in i\R $ \\ \hline
 $ z=t,y<t<-x        $&$ g\in i\R    $&$ dh(\dot{z})\in i\R $ \\ \hline
 $ z=it,\ta<t<\cta   $&$ |g|\equiv 1 $&$ dh(\dot{z})\in i\R $ \\ \hline
 $ z=it,-\cta<t<-\ta $&$ |g|\equiv 1 $&$ dh(\dot{z})\in i\R $ \\ \hline
\end{tabular}
\]

From this table, it is immediate to verify that $dh\cdot dg/g$ is real on the expected planar geodesics of $M$, and pure imaginary on the expected straight lines of $M$.
\\
\\
{\bf 6. Solution of the period problems}
\\ 

Let us consider Figure 6. It reproduces Figure 2(a) and its image under $z$ with some special paths indicated there.
\\
\input epsf
\begin{figure} [ht]
\centerline{
\epsfxsize 13cm
\epsfbox{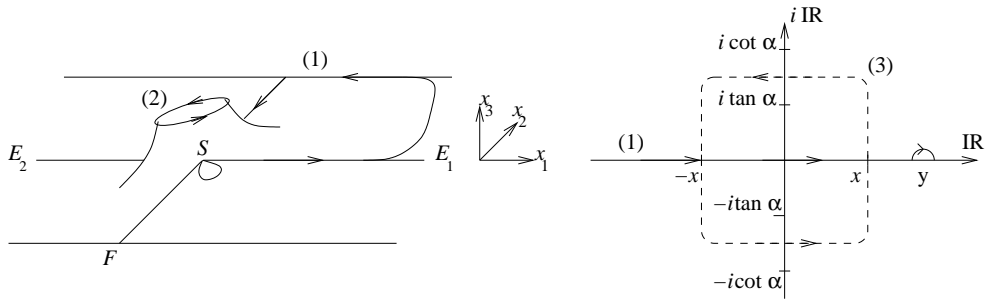}}
\caption{The fundamental piece of $M$ and its image under $z$.}
\end{figure}

Around the punctures of $\ovl{M}$, namely $w^{-1}(\{\pm\cta\})$, we consider small curves given by $w(t)=\pm\cta+\de e^{it}$, where $\de$ is a positive real and $t$ varies in the interval $[0,4\pi]$ (we recall that $w$ takes the values $\pm\cot\af$ with multiplicity 2). All these curves are homotopically equivalent for sufficiently small values of $\de$. Therefore, by letting $\de\to 0$ an immediate calculation leads to $Re\oint(\phi_1,\phi_3)=0$, and up to a minus sign, $Re\oint\phi_2=Res(\phi_2)|_{w=\cta}$, where
\BE
   Res(\phi_2)|_{w=\cta}=\frac{2\pi}{\sqrt{\ctaa-\taa}}\cdot\sqrt{\frac{y+x}{y-x}}.
\EE 

Based on this analysis and (7), it is clear that $Re\int_{(1)}\phi_3=0$. Moreover, (7) also gives us $Re\int_{(2)}\phi_3=0$. The curve (2) from Figure 6 is homotopically equivalent to the sum of (1) with its image under the maps $(g,z)\to(\ovl{g},\ovl{z})$ and $(g,z)\to(-\ovl{g},\ovl{z})$, composed in this order (see Table 4). Actually, this composition corresponds to the rotation $\rho$, explained at the beginning of Section 3. Since $\int_{\rho\circ(1)}(\phi_1,\phi_2)=-\int_{(1)}(\phi_1,\phi_2)$, then $Re\int_{(2)}(\phi_1,\phi_2)=0$. It remains to prove that
\BE
   Re\int_{(1)}(\phi_1,\phi_2)=0.
\EE    

In Figure 6, the curve (3) is symmetric with respect to the geodesic (2). Because of that, the only non-zero component of the period vector $Re\int_{(3)}\phi_{1,2,3}$ must be the third one. Moreover, $Re\int_{(3)}\phi_3\ne 0$ because (7) implies that $dh$ is real and never vanishes on the dashed lines of Figure 3(b). In fact, this component provides the vertical period of $M$, suggested by Figure 1. The horizontal period is given by (8).
\\

We have just reduced the period problems to the proof of (9). For this purpose, we shall apply the {\it limit-method} cited at the introduction. Let us show that the Weierstra\ss \ data (3) and (6) converge to the Weierstra\ss \ data from Callahan-Hoffman-Meeks's surface of genus 3 (see Figure 7).
\\
\begin{figure}
 \centering
 \includegraphics[scale=0.5]{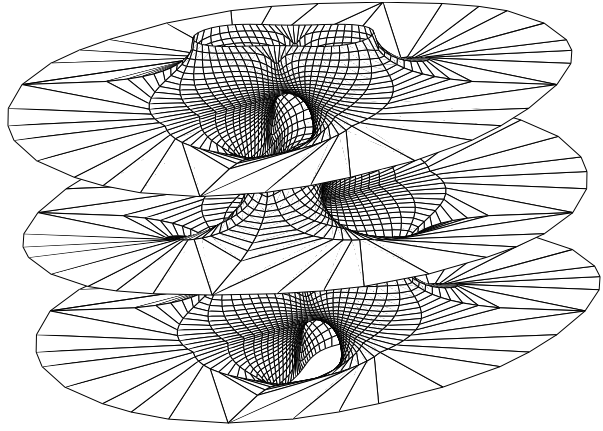}
 \caption{The Callahan-Hoffman-Meeks's surface of genus 3.}
\end{figure}

\begin{figure}[ht]
 \centering
 \includegraphics[scale=0.85]{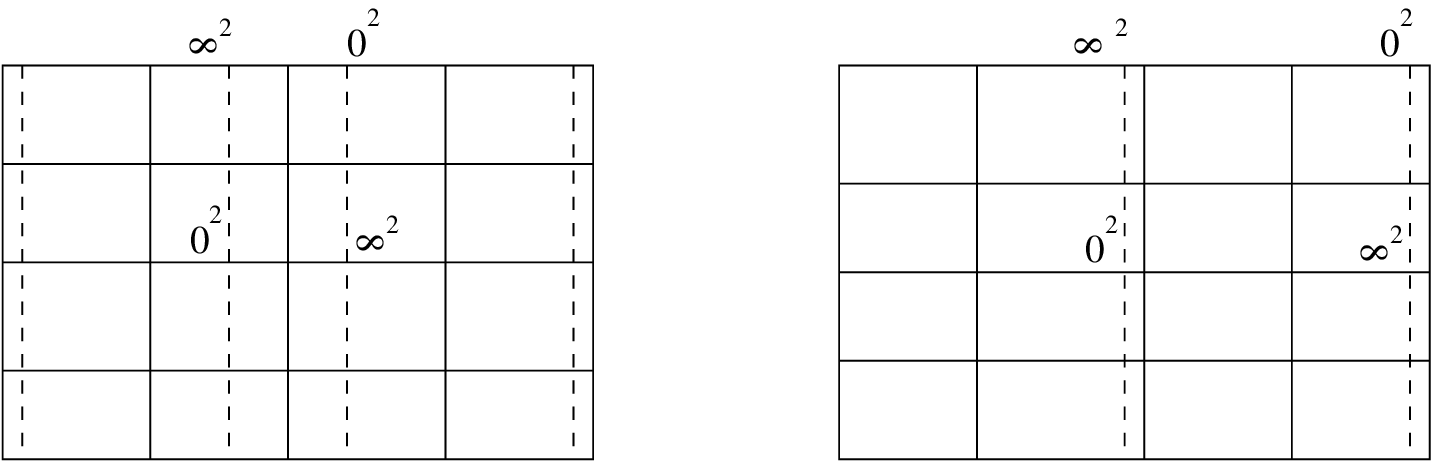}
 \caption{(a) Divisor of $(\frac{x+z}{x-z})^2$; (b) divisor of $(\frac{w-\ta}{w+\ta})(\frac{\cta+w}{\cta-w})$. }
\end{figure}

Consider $K$ a compact subset of $T\setminus w^{-1}(\{\pm\cta\})$. From Figure 8 one sees that $(\frac{x+z}{x-z})^2$ converges uniformly in $K$ to $(\frac{w-\ta}{w+\ta})(\frac{\cta+w}{\cta-w})$ when both $x$ and $y$ approach 1. Thus, from (3) it follows that 
\BE
   g^4=\biggl(\frac{w-\ta}{\ta+w}\biggl)\biggl(\frac{\cta+w}{\cta-w}\biggl)^3.
\EE

By comparing (7) and (10) with [CHM, p. 502] one sees that our surfaces coincide in the limit. This is the first step to solve (9).
\\
REMARK 6.1: If $\phi_{1,2,3}$ are the Weierstra\ss \ data for $(g,dh)$, we call $\widetilde{\phi}_{1,2,3}$ the ones for $(\widetilde{g}, dh):=(e^{-i \pi/4}g,dh)$. In the case of [CHM] one automatically has $Re\int_{(1)}\widetilde{\phi}_1=0$ due to the additional symmetries. 
\\
\begin{figure}[ht]
 \centering
 \includegraphics[scale=0.55]{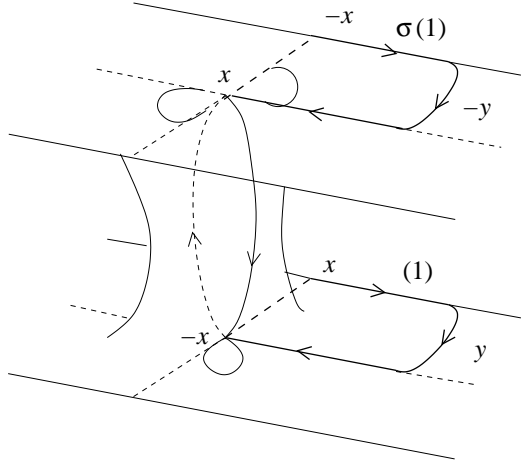}
 \caption{Description of (1) and its image under $\sigma.$}
\end{figure}

Let $\sigma$ be the involution given by $(g,z)\to(-1/\ovl{g},-\ovl{z})$. From (7) and recalling that $w(-\ovl{z})=-\ovl{w}(z)$ we have $dh\to-\ovl{dh}$. Therefore 
$$
   Re \int_{\sigma(1)}\frac{dh}{g}=Re\int_{(1)}\sigma^*\biggl(\frac{dh}{g}\biggl)=Re\int_{(1)}\ovl{gdh}.
$$

Figure (9) suggests that $(1)$ is homotopic to the concatenation of $\beta$ with $\sigma(1)$, where $\beta$ represents the vertical loop and $\sigma(1)$ comes from the involution $\sigma$ applied to $(1)$. This fact can be verified in the complex plane. Hence
\begin{eqnarray*}
  Re\int_{(1)}\phi_1 
  &=& Re\int_{(1)}\frac{dh}{g}-Re\int_{(1)}g \ dh           \\
  &=& Re\int_{(1)}\frac{dh}{g}-Re\int_{\beta\cup\sigma(1)}g \ dh \\ 
  &=& Re\int_{(1)}\frac{dh}{g}-Re\int_{\sigma(1)}g \ dh
      -Re\int_{\beta}g \ dh 
\end{eqnarray*}
\BE
  \hspace{1.7cm}=-Re\int_{\beta}g \ dh, \ {\rm since\eh (1)\eh is\eh a\eh real\eh curve.}
\EE

We split the vertical loop as $\beta:=\beta^+\cup\beta^-$, where $\beta^+$ is the ascending curve from $-x$ to $x$ and $\beta^-$ the path from $x$ to $-x$. The rotation $\rho$, introduced in Section 3, corresponds to $(g,z)\to(-g,z)$, whence $g\to-g$ and $dh\to dh$ under its action. Therefore, 
\BE
   \int_\beta g \ dh=
   \int_{\beta^+\cup\beta^-}g \ dh=2\int_{\beta^+}g \ dh=2\int_{\beta^-}g \ dh.
\EE
\ \\
REMARK 6.2: Figure 10 represents the image under $g$ of a fundamental domain $\D$, namely a smallest subset of $X(M)$ that fully generates it by isometries of $\R^3$. The left image corresponds to $g|_\D$ referring to the ``front piece'', which contains $\beta^-$. The right image concerns the ``back piece'', which contains $\beta^+$.
\\
\begin{figure}
 \centering
 \includegraphics[scale=0.7]{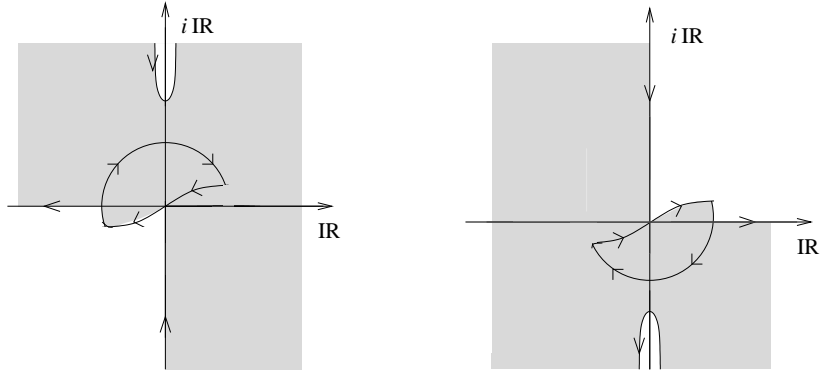}
 \caption{The image under $g$ of a fundamental domain.}
\end{figure}

We consider the following cases:
\\
\\
\underline{Case I} \ $x=1\lessim y$.

For $\beta^+$ we have $z(t)=-e^{it}$ with $0\le t \le\pi$. Hence $dh$ is given by
\begin{eqnarray*}
   dh=\frac{-idz/z}{[z^2+z^{-2}+\taa+\ctaa]^{1/2}}
     =\frac{dt}{[2\cos 2t+\taa+\ctaa]^{1/2}},
\end{eqnarray*}
and $\ta>w(t)>-\ta$. Moreover, $\frac{1+z}{1-z}\big|_{\beta^+}=\frac{-i\sin t}{1+\cos t}$, whence $\frac{\cta+w}{\cta-w}$ and consequently $g^2(t)$ vary according to Figure 11. 
\\

From Remarks 6.1 and 6.2, we notice that the curve in Figure 11(b), rotated by $-\pi/2$, has a branch of square root indicated in Figure 12(a). Therefore $Re\int_{(1)}\widetilde{g}dh<0$.
\\
\begin{figure}
 \centering
 \includegraphics[scale=0.6]{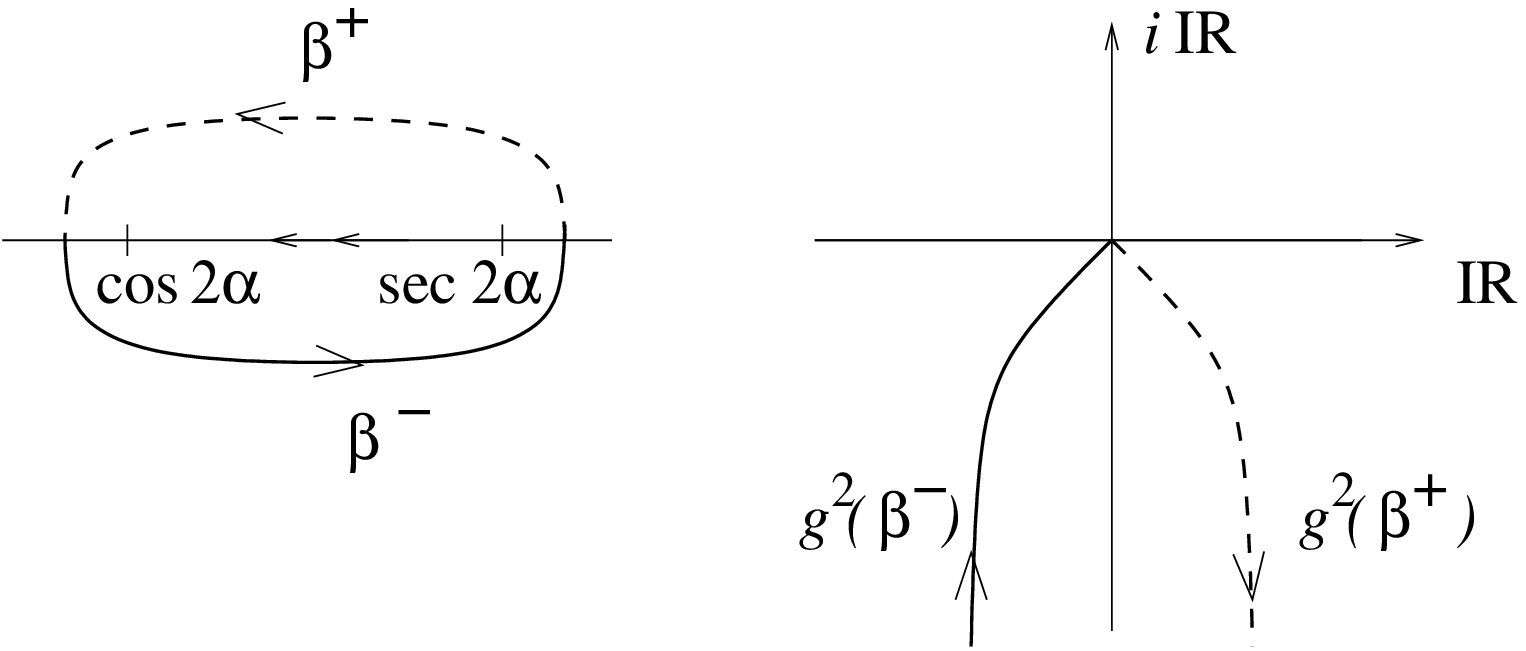}
 \caption{(a) Image of $\frac{\cta+w}{\cta-w}$; (b) image of $g^2$.}
\end{figure}
\begin{figure}
 \centering
 \includegraphics[scale=0.6]{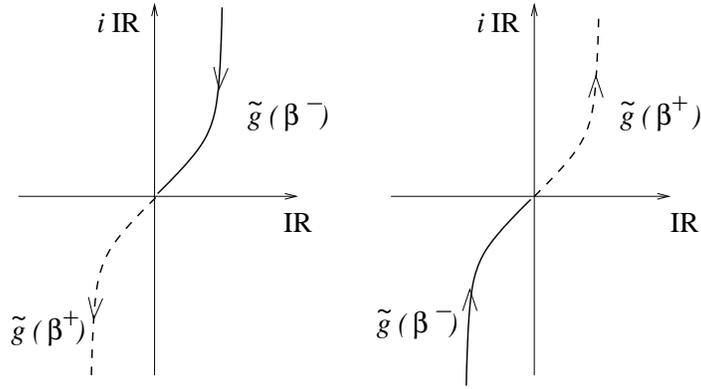}
 \caption{(a) Branch of square root of $\widetilde{g}$ in Case I; (b) in Case II.}
\end{figure}
\begin{figure}
 \centering
 \includegraphics[scale=0.6]{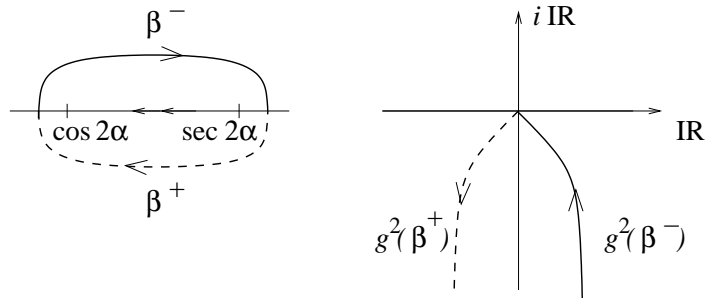}
 \caption{(a) Image of $\beta$ in Case II; (b) image of $g^2$.}
\end{figure}
\ \\
\underline{Case II} \ $x=1\gtrsim y$.

Now $\beta^+$ is still given by $z(t)=-e^{it}$ with $0\le t\le\pi$ and $\frac{1+z}{1-z}=\frac{-i\sin t}{1+\cos t}$. However, since $x\gtrsim y$ we have that $\frac{\cta+w}{\cta-w}$ and consequently $g^2(t)$ vary according to Figure 13. Again from Remarks 6.1 and 6.2, the branch of square root for $\widetilde{g}^2(t)$ is now indicated in Figure 12(b). Thus $Re\int_{(1)}\widetilde{g}dh>0$.
\\

In Figure 14 we indicated the behaviour of a doubly periodic Scherk-Costa surface for the above cases and $x=y=1$.
\\

Therefore, $Re\int_{(1)}\widetilde{g}dh=0$ for some values of $(x,y)$ in a neighbourhood of $(1,1)$. At this limit-point, the function $Re\int _{(1)}\widetilde{\phi}_2$ depends only on the parameter $\af\in(0,\pi /4)$. This is the one-dimensional period problem for [CHM]. According to [MR], it has only {\it one} zero that we call $\af^*$. 
\\
\begin{figure}
 \centering
 \includegraphics[scale=0.55]{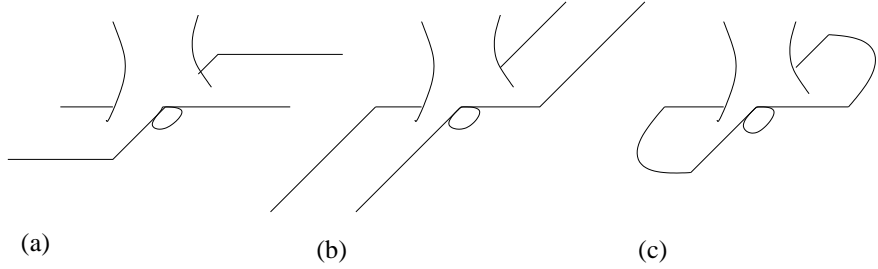}
 \caption{(a) $x<y,$ (b) $x>y,$ (c) $x=y.$}
\end{figure}

We can illustrate this fact by taking a vertical axis $\nu$ and plotting a graph on the plane $O\af\nu$, which crosses the horizontal axis at $\af^*$. Back to our surfaces, if the extra parameters $(x,y)$ were restricted to a curve $(x(\kappa),y(\kappa))$ with an extreme at $(1,1)$, then we could visualize $\kappa$ as a third axis to $O\af\nu$. In this way, both $Re\int _{(1)}\widetilde{\phi}_{1,2}$ turn out to be dependent on two variables, namely $(\af,\kappa)$, and their graphs are surfaces like Figure 15 suggests. Notice that we cannot provide numeric pictures of this fact, since our analyses include limit-values. They typically make unreliable any computational image. This is the second step to solve (9). 
\\
\begin{figure}[ht]
 \centering
 \includegraphics[scale=1]{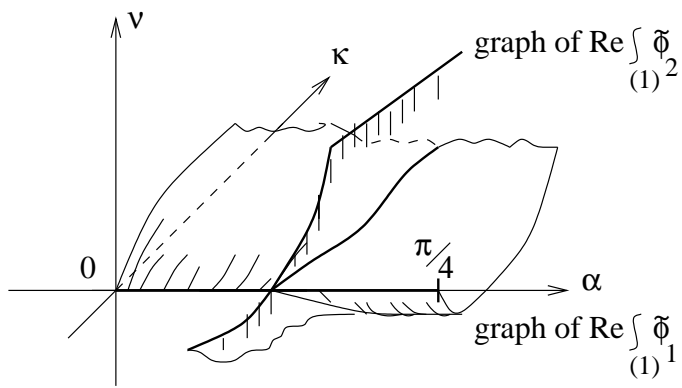}
 \caption{Periods on curve $(1)$.}
\end{figure}

Let us take, for instance, $x$ and $y$ as functions of $(\af,\kappa)$ given by $x=1+\kappa\af$, $y=1+\kappa\af+2\kappa(\frac{4\af^*+\pi}{8}-\af)$. Hence $\kappa\gtrsim 0$ implies $x=1<y$ for $\af\in(0,\frac{4\af^*+\pi}{8}-\af)$ and $x=1>y$ for $\af\in(\frac{4\af^*+\pi}{8}-\af,\frac{\pi}{4})$. Namely, $Re\int_{(1)}\widetilde{\phi}_1<0$ in the first interval and $Re\int_{(1)}\widetilde{\phi}_1>0$ in the second. We could extend $y(\af,\kappa)$ to $y=1+\kappa\af+2\kappa(\frac{4\af^*+\pi}{8}\cdot s-\af)$, $1\ge s>0$. Consequently, there exists a curve $(\af(t),\kappa(t))$ for which $Re\int_{(1)}\widetilde{\phi}_{1,2}=0$. Moreover, along this curve we have $x\ne y$ as explained next.
\\

If $x=y\ne 1$, we assert that the period $Re\int_{(1)}\widetilde{\phi}_1$ is {\it non-zero} in the $x_1$-direction. This is because one gets a CHM-surface with ``torsion'', as illustrated in Figure 16(b). Without torsion, on $\beta$ one has real $dh$ and $\widetilde{g}=-i|\widetilde{g}|$, whence $Re\int_{\beta}\widetilde{g}dh=0$. For $x=y\ne 1$, however, we may still set $dh|_{\beta^+}$ to be real and positive, but then $\widetilde{g}|_{\beta^+}$ gets a never-vanishing real part. Therefore, from (11) and (12) one has $Re\int_{(1)}\widetilde{\phi}_1dh\ne 0$. This third step finally proves (9) and concludes the present section.
\\
\begin{figure}
 \centering
 \includegraphics[scale=0.8]{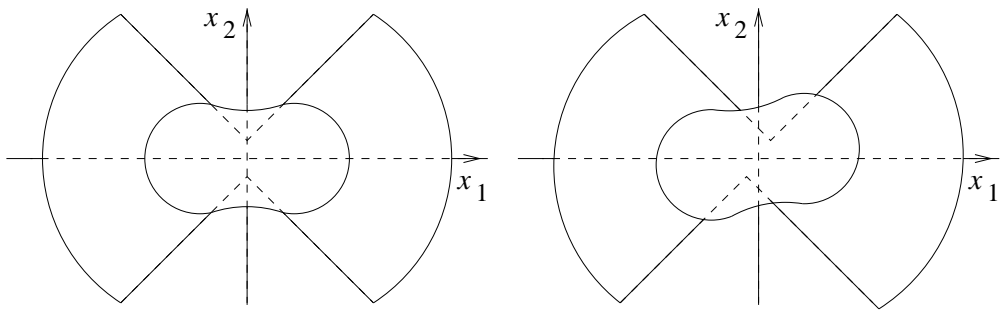}
 \caption{Symmetries for $x=y=1$ (CHM) and $x=y\ne 1$ (CHM with torsion).}
\end{figure}
\ \\
{\bf 7. Embeddedness of the fundamental piece}
\\

This chapter is strongly based in the ideas of [MRB] and [RB4]. We begin with by identifying a fundamental domain $\D$ of $X(M)$ in Figure 17.   
\\
\begin{figure}[ht]
 \centering
 \includegraphics[scale=0.8]{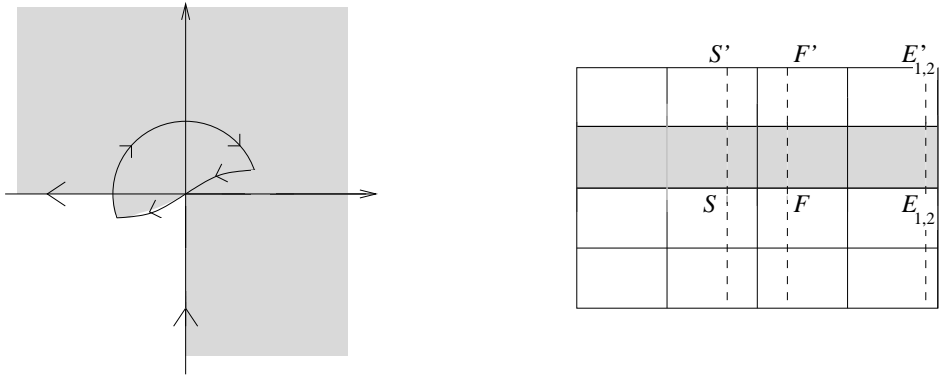}
 \caption{(a) The set $g(\D)$; (b) the image of the fundamental domain $\D$ under $\rho$.}
\end{figure}

In the previous section we proved the existence of a curve $(x(t),y(t))$, $0\le t<1$, along which (9) holds. Moreover, $\Lim{t\to0}{\af(t)}=\af^*$. Let us fix $t\in(0,1)$ and consider the minimal immersion $X_t:\D\setminus\{E_{1,2}\}\to\R^3$, defined by (3) and (6). Each branch of square root of $g$ takes any $q\in\D\setminus\{E_{1,2}\}$ to a pair of points in $\R^3$, say $X_t(q)^+$ and $X_t(q)^-$. If $X_t(0)$ is the origin, then each point is the image of the other by a $180^\circ$ rotation about $Ox_3$.
\\

We consider a fundamental piece $P$ of $M$. Let $P^-$ be the image of $\D\setminus\{E_{1,2}\}$ in $\R^3$ under $X_t$, and $P^+$ the image of $P^-$ in $\R^3$ under a $180^\circ$ rotation around $X_t([-x,x])$. Thus $P=P^+\cup P^-$. 
\\

Let $\K$ be a subset of $\mathcal{D}$ such that $\D\setminus\K=V_E$, where $V_E$ is a connected neighbourhood of $E_{1,2}$. From (7) and (10) we see that $(g,dh)$ converge uniformly to the Weierstra\ss \ data of the embedded CHM-surface. Let us denote this minimal embedding by $X_0$. When $t\to 0$, $X_t|_{V_E}$ approaches a planar end for sufficiently small $V_E$. For $t$ close to zero, the projection of $X_t|_{\partial V_E}$ onto $x_3=0$ consists of two curves $C^\pm$ which determine two simply connected open regions $R^+$ and $R^-$. Since $g(V_E)$ is contained in a half-sphere, then $(x_1,x_2)|_{V_E}$ is an immersion onto $R^{\pm}$ because $x_2$ is bounded for any fixed $t\in(0,1)$. Since $\partial R^\pm$ are the monotone curves $C^\pm$, then $X_T|_{\partial V_E}$ is a graph of $x_3$ as a function of $(x_1,x_2)$.
\\
\begin{figure}[ht]
 \centering
 \includegraphics[scale=0.6]{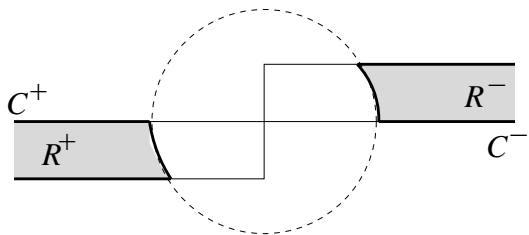}
 \caption{Regions $R^\pm$ and curves $C^\pm$.}
\end{figure}

We observe that $X_0|_\K$ is a compact embedded minimal surface in $\R^3$. Since its boundary does not have self-intersections, then $X_t|_\K$ is still embedded for sufficiently small $t$. Moreover, $X_t|_\K$ does not intercept $X_t|_{V_E}$, otherwise there would be a ball in $\R^3$ containing the whole boundary of $X_t|_\K$ but not all the rest of it. This is impossible according to the {\it maximum principle}. Hence, the pieces $X_t|_\K$ and $X_t|_{V_E}$ make together a minimal embedding $X_t:\D\setminus\{E_{1,2}\}\to\R^3$, for $t$ sufficiently close to zero.
\\

Again by the {\it maximum principle}, we may extend this conclusion for all $t\in (0,1)$. Therefore, $P^+$ is embedded in $\R^3$, and since $P^-$ is its image under a $180^\circ$ rotation about the segment of $P^+$, the whole piece $P$ will not have self-intersections. Since the immersion is proper, then $P$ is embedded in $\R^3$.
\\

Now $P\subset\R^3/\G$, where $\G$ is the group of $\R^3$ generated by $(x_1,x_2,x_3)\to(x_1,x_2,-x_3+2Re\int_{\beta^+}dh)$ and $(x_1,x_2,x_3)\to(x_1,x_2+Re\oint\phi_2,x_3)$. In the horizontal faces of $\partial(\R^3/\G)$ we have the reflection curves of $P$. In the vertical faces we have the straight lines of $P$. By applying $\G$ to $P$ one generates $M$, which is then complete, doubly periodic and embedded in $\R^3$.
\ \\
\ \\
{\bf References}
\ \\
\ \\
$[$BRB$]$-F. Baginski and V. Ramos Batista: Solving period problems for minimal surfaces with the support function. In manuscript (2007); home page http://www.ufabc.edu.br/pgmatematica/docentes.html
\\
$[$CHM$]$-M. Callahan, D. Hoffman and W.H. Meeks: Embedded minimal surfaces with an infinite number of ends. Inventiones Math. {\bf 96} (1989), 459-505.
\\
$[$H$]$-G. Hart: Where are nature's missing structures? Nature Materials {\bf 6} (2007), 941--945.
\\
$[$HKW$]$-D. Hoffman, H. Karcher and F. Wei: The genus one helicoid and the minimal surfaces that led to its discovery. Global Analysis and Modern Mathematics. Publish or Perish Press (1993), 119--170.  
\\
$[$JS$]$-H. Jenkins and J. Serrin: Variational problems of minimal surface type. II. Boundary value problems for the minimal surface equation. Arch. Rational Mech. Anal. {\bf 21} (1966), 321--342.
\\
$[$K2$]$-H. Karcher: Construction of minimal surfaces, Surveys in Geometry, University of Tokyo (1989), 1--96, and Lecture Notes {\bf 12}, SFB256, Bonn (1989).
\\
$[$K1$]$-H. Karcher: Embedded minimal surfaces derived from Sckerk's examples. Manuscr. Math. {\bf 62} (1988), 83--114.
\\
$[$Kp$]$-N. Kapouleas: Complete embedded minimal surfaces of finite total curvature. J. Differential Geom. {\bf 47} (1997), 95--169. 
\\
$[$LMa$]$-F. L\'opez and F. Mart\ih n: Complete minimal surfaces in $\R^3$. Publicacions Matematiques {\bf 43} (1999), 341--449.
\\
$[$LM$]$-E. Lord and A. Mackay: Periodic minimal surfaces of cubic symmetry. Current Science {\bf 85} (2003), 346--362. 
\\
$[$L$]$-K. L\"ubeck: M\'etodo-limite para solu\c c\~ao de problemas de per\ih odos em superf\ih cies m\ih nimas. Doctoral Thesis, University of Campinas (2007).
\\
$[$LRB$]$-K. L\"ubeck and V. Ramos Batista: A limit-method for solving period problems on minimal surfaces. In manuscript (2007); home page 

\hspace{-0.6cm}http://www.ufabc.edu.br/pgmatematica/docentes.html
\\
$[$MRB$]$-F. Mart\ih n and V. Ramos Batista: The embedded singly periodic Scherk-Costa surfaces. Math. Ann. {\bf 336} (2006), 1, 155--189. 
\\
$[$MR$]$-F. Mart\ih n and D. Rodr\ih guez: A characterization of the periodic Callahan-Hoffman-Meeks surfaces in terms of their symmetries. Duke Math. J. {\bf 89} (1997), 445--463.
\\
$[$N$]$-J. Nitsche: Lectures on minimal surfaces. Cambridge University Press, Cambridge (1989).
\\
$[$O$]$-R. Osserman: A survey of minimal surfaces, Dover, New York, 2nd ed (1986).
\\
$[$PRT$]$-J. Perez, M. Rodr\ih iguez and M. Traizet: The classification of doubly periodic minimal tori with parallel ends. J. Differential Geom. {\bf 69} (2005), 523-­577.
\\
$[$RB4$]$-V. Ramos Batista: Singly periodic Costa surfaces. J. London Math. Soc. {\bf 72} (2005), 2, 478--496. 
\\
$[$RB3$]$-V. Ramos Batista: Noncongruent minimal surfaces with the same symmetries and conformal structure. Tohoku Math. J. {\bf 56} (2004), 237--254.
\\
$[$RB2$]$-V. Ramos Batista: A family of triply periodic Costa surfaces, Pacific J. Math. {\bf 212} (2003), 347--370.
\\
$[$RB0$]$-V. Ramos Batista: Construction of new complete minimal surfaces in $\R^3$ based on the Costa surface. Doctoral Thesis, University of Bonn (2000).
\\
$[$T3$]$-M. Traizet: An embedded minimal surface with no symmetries. J. Differential Geom. {\bf 60} (2002), 103--153.
\\
$[$T2$]$-M. Traizet: Adding handles to Riemann's minimal surfaces. J. Inst. Math. Jussieu {\bf 1} (2002), 145--174.
\\
$[$T1$]$-M. Traizet: Construction de surfaces minimales en recollant des surfaces de Scherk. Ann. Inst. Fourier {\bf 46} (1996), 1385--1442.
\\
$[$Web2$]$-M. Weber: A Teichm\"uller theoretical construction of high genus singly periodic minimal surfaces invariant under a translation. Manuscripta Math. {\bf 101} (2000), 125--142.
\\
$[$Web1$]$-M. Weber: The genus one helicoid is embedded. Habilitation Thesis, Bonn 2000.
\\
$[$W$]$-F. Wei: Some existence and uniqueness theorems for doubly periodic minimal surfaces. Invent. Math. {\bf 109} (1992), 113--136.
\\ \\
L\"ubeck, Kelly\\
Universidade Estadual do Oeste do Paran\'a\\
av. Tarquinio Joslin dos Santos, 1300\\
P.O.Box 961\\
85870-650 Foz do Igua\c cu - PR, Brazil\\
{\tt klubeck@unioeste.br}\\
{\tt http://www.foz.unioeste.br/$\sim$lem/Docentes\%20Kelly.htm}\\ \\
Ramos Batista, Val\'erio\\
Universidade Federal do ABC\\
r. Catequese 242, 3rd floor\\
09090-400 Santo Andr\'e - SP, Brazil\\
{\tt valerio.batista@ufabc.edu.br}\\ 
{\tt http://www.ufabc.edu.br/pgmatematica}

\end{document}